
\magnification1200
\input amstex.tex
\documentstyle{amsppt}

\hsize=12.5cm
\vsize=18cm
\hoffset=1cm
\voffset=2cm

\footline={\hss{\vbox to 2cm{\vfil\hbox{\rm\folio}}}\hss}
\nopagenumbers
\def\DJ{\leavevmode\setbox0=\hbox{D}\kern0pt\rlap
{\kern.04em\raise.188\ht0\hbox{-}}D}

\baselineskip=13pt
\def\hf{{\textstyle{1\over2}}}
\def\a{\alpha}
\def\d{{\,\roman d}}
\def\e{\varepsilon}\def\E{{\roman e}}
\def\f{\varphi}
\def\G{\Gamma}

\def\={\;=\;}
\def\zx{\zeta(\hf+ix)}
\def\zt{\zeta(\hf+it)}
\def\zu{\zeta(\hf+iu)}
\def\D{\Delta}
\def\E{{\roman e}}

\def\z{\zeta}

\def\hf{{\textstyle{1\over2}}}

\def\f{\varphi}

\def\le{\leqslant} \def\ge{\geqslant}
\font\tenmsb=msbm10
\font\sevenmsb=msbm7
\font\fivemsb=msbm5
\newfam\msbfam
\textfont\msbfam=\tenmsb
\scriptfont\msbfam=\sevenmsb
\scriptscriptfont\msbfam=\fivemsb
\def\Bbb#1{{\fam\msbfam #1}}

\def \NN {\Bbb N}

\font\ff=cmr8

\baselineskip=13pt

\font\teneufm=eufm10
\font\seveneufm=eufm7
\font\fiveeufm=eufm5
\newfam\eufmfam
\textfont\eufmfam=\teneufm
\scriptfont\eufmfam=\seveneufm
\scriptscriptfont\eufmfam=\fiveeufm
\def\mathfrak#1{{\fam\eufmfam\relax#1}}

\font\tenmsb=msbm10
\font\sevenmsb=msbm7
\font\fivemsb=msbm5
\newfam\msbfam
     \textfont\msbfam=\tenmsb
      \scriptfont\msbfam=\sevenmsb
      \scriptscriptfont\msbfam=\fivemsb
\def\Bbb#1{{\fam\msbfam #1}}

\def \NN {\Bbb N}

  \def\rightheadline{{\hfil{\ff
  On zeta estimates in short intervals}\hfil\tenrm\folio}}

  \def\leftheadline{{\tenrm\folio\hfil{\ff
   Aleksandar Ivi\'c }\hfil}}
  \def\emptyheadline{\hfil}
  \headline{\ifnum\pageno=1 \emptyheadline\else
  \ifodd\pageno \rightheadline \else \leftheadline\fi\fi}

\topmatter
\title
ON SOME MEAN VALUE RESULTS FOR THE ZETA-FUNCTION IN SHORT INTERVALS
\endtitle
\author   Aleksandar Ivi\'c  \endauthor

\bigskip
\address
Aleksandar Ivi\'c, Katedra Matematike RGF-a
Universiteta u Beogradu, \DJ u\v sina 7, 11000 Beograd, Serbia
\endaddress
\email {\tt
ivic\@rgf.bg.ac.rs,  aivic\@matf.bg.ac.rs} \endemail
\dedicatory
\enddedicatory
\abstract
{Let $\D(x)$ denote the error term in the Dirichlet
divisor problem, and let $E(T)$ denote the error term in the asymptotic
formula for the mean square of $|\zt|$. If
$E^*(t) := E(t) - 2\pi\D^*(t/(2\pi))$ with $\D^*(x) =
 -\D(x)  + 2\D(2x) - \hf\D(4x)$ and
 $\int_0^T E^*(t)\d t = \frac{3}{4}\pi T + R(T)$, then we obtain
a number of results involving the moments of $|\zt|$
in short intervals, by connecting
them to the moments of $E^*(T)$ and $R(T)$ in short intervals. Upper bpunds
and asymptotic formulas for integrals
of the form
$$
\int_T^{2T}\left(\int_{t-H}^{t+H}|\zu|^2\d u\right)^k\d t\qquad(k\in\NN, 1 \ll H \le T)
$$
are also treated.}
\endabstract
\endtopmatter

\document
\bigskip
\head
1. Introduction
\endhead
\bigskip

As usual, let
$$
\D(x) \;:=\; \sum_{n\le x}d(n) - x(\log x + 2\gamma - 1)
\leqno(1.1)
$$
denote the error term in the classical Dirichlet divisor problem. Also let
$$
E(T) \;:=\;\int_0^T|\zt|^2\d t - T\left(\log\bigl({T\over2\pi}\bigr) + 2\gamma - 1
\right)\leqno(1.2)
$$
denote the error term in the mean square formula for $|\zt|$.
Here $d(n)$ is the number of divisors of
$n$, $\z(s)$ is the Riemann zeta-function, and $ \gamma = -\G'(1) = 0,577215\ldots\,$
is Euler's constant. In view of F.V. Atkinson's classical explicit formula
for $E(T)$ (see [1],  [4, Chapter 15] and [5, Chapter 2]) it was known long ago that
there are analogies between $\D(x)$ and $E(T)$. However, in this context it seems
that instead of the error-term function $\D(x)$ it
is more exact to work with the modified
function $\D^*(x)$ (see  M. Jutila [14], [15] and T. Meurman [17]), where
$$
\eqalign{
\D^*(x) :&= -\D(x)  + 2\D(2x) - \hf\D(4x)\cr&
= \hf\sum_{n\le4x}(-1)^nd(n) - x(\log x + 2\gamma - 1),\cr}
\leqno(1.3)
$$
since it turns out that $\D^*(x)$ is a better analogue of $E(T)$ than $\D(x)$.
Namely, M. Jutila (op. cit.) investigated both the
local and global behaviour of the difference
$$
E^*(t) \;:=\; E(t) - 2\pi\D^*\bigl({t\over2\pi}\bigr),
$$
and in particular in [14] he proved that
$$
\int_T^{T+H}(E^*(t))^2\d t \;\ll_\e\; HT^{1/3}\log^3T+ T^{1+\e}\quad(1\le H\le T).\leqno(1.4)
$$
Here and later $\e$ denotes positive constants which are arbitrarily
small, but are not necessarily the same ones at each occurrence,
while $a \ll_\e b$ (same as $a = O_\e(b))$ means that
the $\ll$--constant depends on $\e$. The significance of (1.4) is that, in view of
(see e.g., [4])
$$
\eqalign{
\int_0^T(\D^*(t))^2\d t &\;\sim\; AT^{3/2},\cr
 \int_0^T E^2(t)\d t &\;\sim \;BT^{3/2}\qquad(A,B >0,\;
T\to\infty),\cr}
$$
it transpires that $E^*(t)$ is in the mean square sense of a lower order of magnitude than
either $\D^*(t)$ or $E(t)$.

Later works provided more results on the mean values of $E^*(T)$. Thus
in [9] the author sharpened (1.4) (in the case when $H=T$) to the asymptotic formula
$$
\int_0^T (E^*(t))^2\d t \;=\; T^{4/3}P_3(\log T) + O_\e(T^{7/6+\e}),\leqno(1.5)
$$
where $P_3(y)$ is a polynomial of degree three in $y$ with
positive leading coefficient, and all the coefficients may be evaluated
explicitly.  This, in particular, shows that (1.4) may be complemented with the
lower bound
$$
\int_T^{T+H}(E^*(t))^2\d t \;\gg\; HT^{1/3}\log^3T\qquad(T^{5/6+\e} \le H \le T),\leqno(1.6)
$$
which is implied by (1.5).
It seems likely that the error term in (1.5) is $O_\e(T^{1+\e})$, but this seems difficult
to prove. In [12] the author showed that (1.6) remains true for $T^{2/3+\e} \le H \le T$.

\medskip
In what concerns higher moments of $E^*(T)$ the author proved ([6,  Part 4])
$$
\int_0^T |E^*(t)|^3\d t \;\ll_\e\; T^{3/2+\e},\leqno(1.7)
$$
in [6, Part 2] that
$$
\int_0^T |E^*(t)|^5\d t \;\ll_\e\; T^{2+\e},\leqno(1.8)
$$
so that by the Cauchy-Schwarz for integrals (1.7) and (1.8) yield
$$
\int_0^T (E^*(t))^4\d t \;\ll_\e\; T^{7/4+\e}.\leqno(1.9)
$$

\medskip
In part [6, Part 3]  the error-term function $R(T)$ was  introduced by
 the relation
$$
\int_0^T E^*(t)\d t \= {3\pi\over 4}T + R(T).\leqno(1.10)
$$
It was shown, by using an estimate for two-dimensional exponential sums, that
$$
R(T) = O_\e(T^{593/912+\e}), \quad {593\over912} = 0,6502129\ldots\,.
\leqno(1.11)
$$
In the same paper it was  also proved that

\bigskip

$$
\int_0^TR^2(t)\d t \= T^2p_3(\log T) + O_\e(T^{11/6+\e}),\leqno(1.12)
$$
where $p_3(y)$ is a cubic polynomial in $y$ with positive leading
coefficient, whose all coefficients may be explicitly evaluated,
and
$$
\int_0^TR^4(t)\d t \;\ll_\e\; T^{3+\e}.
\leqno(1.13)
$$

\medskip
The asymptotic formula (1.12) bears resemblance to (1.5), and it is proved
by a similar technique. The exponents in the error terms are, in both cases,
less than the exponent of $T$ in the main term by 1/6. From (1.5) one obtains that $E^*(T)
= \Omega\bigl(T^{1/6}(\log T)^{3/2}\bigr)$, which shows that $E^*(T)$ cannot be
too small ($f(x) = \Omega(g(x))$ means that $f(x) = o(g(x))$ does not hold
as $x\to\infty$). Likewise, (1.12) yields
$$
R(T) \= \Omega\Bigl(T^{1/2}(\log T)^{3/2}\Bigr).
\leqno(1.14)
$$
It seems plausible that the error term in (1.12) should be $O_\e(T^{5/3+\e})$,
and one may conjecture that
$$
R(T) \= O_\e(T^{1/2+\e})\leqno(1.15)
$$
holds, which is supported by (1.12). In [12] it was proved that,
in the range $T^{2/3+\e} \le H \le T$, we have
$$
\int_T^{T+H}R^2(t)\d t \;\gg \; HT\log^3T,\leqno(1.16)
$$
 and, for $T^\e \le H \le T$,
$$
\int_T^{T+H}R^2(t)\d t \;\ll_\e\; HT\log^3T + T^{5/3+\e}.\leqno(1.17)
$$

\bigskip
\head
2. Statement of results
\endhead
\bigskip

Mean values (or moments) of $|\zt|$ represent  one of the central themes in the theory of $\z(s)$.
There are two monographs dedicated solely to it: the author's [5], and that of K.
Ramachandra [18].
Our results connect bounds for the moments of $|\zt|, E^*(t)$ and $R(t)$ in short intervals.
The meaning of ``short interval'' is that $[T, T+H]$ is such an interval where one can have $H$
much smaller than $T$, namely $H = o(T)$ as $T\to\infty$.
The results are contained in
\bigskip
THEOREM 1. {\it For $k\in \NN$ fixed, $T^{1/3}\le H = H(T)\le T$, we have
$$
\int_T^{T+H}|\zt|^{2k+2}\d t \;\ll_k\; (\log T)^{k+2}\int_{T-H}^{T+2H}|E^*(t)|^{k}\d t. \leqno(2.1)
$$
and }
$$
\int_T^{T+H}|E^*(t)|^{2k}\d t \;\ll_k\; (\log T)^{k+2}\int_{T-H}^{T+2H}|R(t)|^{k}\d t. \leqno(2.2)
$$

\bigskip

THEOREM 2. {\it Let $k\in\NN$ be fixed and $T^\e \le H=H(T)\le T$. If
$$
\int_0^T|E^*(t)|^k\d t \;\ll_{\e,k}\; T^{A(k)+\e}\leqno(2.3)
$$
for some constant $A(k)$, then we must have $A(k) \ge 1 + k/6$, and}
$$
\int_T^{2T}\left(\int_{t-H}^{t+H}|\zu|^2\d u\right)^k\d t \;\ll_{\e,k}\; T^{A(k)+\e} +TH^k(\log T)^{k}.
\leqno(2.4)
$$

\medskip
When $k=1$ or $k=2$ a much more precise result can be obtained for the integral in (2.4). This
is contained in

\medskip
THEOREM 3. {\it
For $T^\e \le H =H(T)\le T$ we have
$$
\int_T^{2T}\Bigl(\int_{t-H}^{t+H}|\zu|^2\d u\Bigr)\d t =
2H\left(T\log \Bigl(\frac{2T}{\E\pi} \Bigr)\right) + O(H^2) + O(T^{3/4}).\leqno(2.5)
$$
and
$$
\int_T^{2T}\left(\int_{t-H}^{t+H}|\zu|^2\d u\right)^2\d t \;\ll \; H^2T(\log T)^2.\leqno(2.6)
$$
For $T^\e \le H = H(T) \le T^{1/2-\e}$ we have the asymptotic formula
$$
\eqalign{&
\int_T^{2T}\left(\int_{t-H}^{t+H}|\zu|^2\d u\right)^2\d t = H^2T(4\log ^2T + e_1\log T + e_0)\cr&
+ HT\sum_{j=0}^3d_j\log^j\Bigl({\sqrt{T}\over 2H}\Bigr)+ O_\e(T^{1/2+\e}H^2)
+ O_\e(T^{1+\e}H^{1/2}).\cr}
\leqno(2.7)
$$
where the $d_j$'s and $e_1, e_0$ are suitable constants ($d_3>0$).}
\bigskip

\medskip
The proofs of Theorem 1, Theorem 2  and Theorem 3 will be given  in Section 3.
In Section 4 we shall provide some corollaries and remarks to these theorems.

\bigskip
\head
3. Proofs of the Theorems
\endhead

\bigskip
In (2.1) of Theorem 1 we have an estimate for the moments of $|\zt|$.
In order to deal with these moments we shall use the standard
large values technique (see e.g.,  [4, Chapter 8]). To transform discrete sums
into sums of integrals one uses the bound
$$
|\zt|^k \ll \log t\int_{t-1}^{t+1}|\zx|^k \d x + 1,
\qquad(k\in \NN\;\;{\roman{fixed}})
\leqno(3.1)
$$
which is Theorem 1.2 of [5] (see also Lemma 7.1 of [4]).
\bigskip
We begin (henceforth let $L = \log T$ for brevity) by noting that, for $T^\e \ll G \le T$,
$$
\eqalign{&
\int_{T-G}^{T+G}|\zt|^2\d t = \int_{-G}^G |\zeta(\hf+iT+iu)|^2\d u \cr&\le
\E \int_{-\infty}^\infty |\zeta(\hf + iT+iu)|^2\E^{-(u/G)^2}\d u\cr&
= \E \int_{-GL}^{GL}|\zeta(\hf + iT+iu)|^2\E^{-(u/G)^2}\d u + O(\E^{-{1\over2} L^2}).
\cr}
$$
In view of (1.2) we further have, on integrating by parts,
$$
\eqalign{&
\int_{-GL}^{GL}|\zeta(\hf + iT+iu)|^2\E^{-(u/G)^2}\d u =
\int_{-GL}^{GL}\E^{-(u/G)^2}\d E(T+u) + O(GL)\cr&
= 2\int_{-GL}^{GL}uG^{-2}\E^{-(u/G)^2}E(T+u)\d u + O(GL).\cr}
$$
By the definition of $E^*(T)$ the last integral becomes
$$
\frac{1}{G^2}\int_{-GL}^{GL}xE^*(T+x)\E^{-(x/G)^2}\d x
+ \frac{2\pi}{G^2}\int_{-GL}^{GL}x\D^*\left(\frac{T+x}{2\pi}\right)\E^{-(x/G)^2}\d x.
$$
To bound the integral containing the $\D^*$ function, we shall use the estimate
$$
\sum_{x<n\le x+h}d(n) \ll h\log x\qquad(x^\e \le h \le x),\leqno(3.2)
$$
which follows from a general result of P. Shiu [19] on multiplicative functions.
Write
$$
\int_{-GL}^{GL}x\D^*\left(\frac{T+x}{2\pi}\right)\E^{-(x/G)^2}\d x
= \int_{-GL}^0\cdots \d x + \int_{0}^{GL}\cdots \d x,\leqno(3.3)
$$
and make the change of variable $y = -x$ in the first integral on the right-hand side.
Then (3.3) becomes
$$
\eqalign{&
\int_0^{GL}y\D^*\left(\frac{T-y}{2\pi}\right)\E^{-(y/G)^2}\d y
+ \int_{0}^{GL}x\D^*\left(\frac{T+x}{2\pi}\right)\E^{-(x/G)^2}\d x\cr&
= \int_{0}^{GL}x\left\{\D^*\left(\frac{T+x}{2\pi}\right)-
\D^*\left(\frac{T-x}{2\pi}\right)\right\}\E^{-(x/G)^2}\d x.\cr}
$$
For $|x| \le T^{\e/3}$ we use the trivial bound (coming from $d(n) \ll_\e n^{\e/3}$)
$$
\D^*\left(\frac{T+x}{2\pi}\right)- \D^*\left(\frac{T-x}{2\pi}\right)
\ll_\e T^{2\e/3},
$$
while for $T^{\e/3} < |x| \le GL$ we use (3.2). This yields
$$
\eqalign{&
\int_{0}^{GL}x\left\{\D^*\left(\frac{T+x}{2\pi}\right)-
\D^*\left(\frac{T-x}{2\pi}\right)\right\}\E^{-(x/G)^2}\d x\cr&
\ll_\e T^{2\e/3}G^2 + G^3L \ll G^3L,\cr}\leqno(3.4)
$$
since $G \gg T^\e$. Therefore  (3.4) furnishes the bound
$$
\frac{2\pi}{G^2}\int_{-GL}^{GL}x\D^*\left(\frac{T+x}{2\pi}\right)\E^{-(x/G)^2}\d x
\;\ll\; GL,
$$
and we obtain the starting point for the proof of (2.1),
which we formulate as

\medskip
LEMMA 1. For $T^\e \le G = G(T)\le T, L = \log T$ we have
$$
\int_{T-G}^{T+G}|\zt|^2\d t \le \frac{2\E }{G^2}\int_{-GL}^{GL}xE^*(T+x)\E^{-(x/G)^2}\d x
+ O(GL).\leqno(3.5)
$$

\medskip
We return to the proof of (2.1) and
suppose now that $\{t_r\}_{r=1}^R$ is a set of points satisfying
$$
T < t_1 < \ldots < t_R \le T+H,\; |\zeta(\hf+it_r)| \ge V,\; |t_r-t_s| \ge 1
\;(r = 1, \ldots, R).\leqno(3.6)
$$
We use (3.1) and group the intervals $[t_r-1, t_r+1]$ into disjoint subintervals of the form
$$
[\tau_s-G, \tau_s+G]\quad(s = 1,2,\ldots,  S\le R,\; G \ll H).
$$
Then by Lemma 1 we obtain (we may suppose that the sum over $s$ below is the largest
of several sums of the same type)
$$
\eqalign{
RV^{2} &\ll L\sum_{s=1}^S\int_{\tau_s-G}^{\tau_s+G}|\zt|^2\d t\cr&
\le 2\E L\sum_{s=1}^S\frac{1}{G^2}\int_{-GL}^{GL}xE^*(\tau_s+x)\E^{-(x/G)^2}\d x,
\cr}\leqno(3.7)
$$
provided that, for some sufficiently small $c>0$, we choose
$$
G \;=\; cV^2/L.\leqno(3.8)
$$
By bounds for $|\zt|$ we obtain $G \ll T^{1/3} \ll H$, and we
choose a representative set of points $\tau_\ell, \ell = 1, \ldots, S'\,(\le S)$
from the set $\{\tau_s\}_{s=1}^S$ such that
the intervals $(\tau_\ell-GL, \,\tau_\ell+GL)$ are disjoint for $\ell = 1,\ldots, S'$.
Therefore it follows by H\"older's  inequality for integrals that
$$
\eqalign{&RV^{2} \ll
L\sum_{\ell=1}^{S'}G^{-2}\int_{-GL}^{GL}|xE^*(\tau_\ell+x)|\E^{-(x/G)^2}\d x\cr&
\ll LG^{-2}\sum_{\ell=1}^{S'}\left(\int\limits_{-GL}^{GL}|E^*(\tau_\ell+x)|^k
\E^{-(x/G)^2}\d x\right)^{\frac{1}{k}}
\left(\int\limits_{-GL}^{GL}|x|^{k\over k-1}\E^{-(x/G)^2}\d x\right)^{1-\frac{1}{k}}\cr&
\ll LG^{-2}\left(\sum_{\ell=1}^{S'}\int_{-GL}^{GL}|E^*(\tau_\ell+x)|^k
\E^{-(x/G)^2}\d x\right)^{\frac{1}{k}}
G^{2-\frac{1}{k}}(S')^{1-\frac{1}{k}}\cr&
\ll G^{-\frac{1}{k}}LS^{1-\frac{1}{k}}\left(\int_{T-H}^{T+2H}|E^*(x)|^k\d x\right)^{\frac{1}{k}}.
\cr&}
$$
Since $S\le R$, in view of (3.8) this gives
$$
\eqalign{
R &\ll V^{-2k}G^{-1}L^{k}\int_{T-H}^{T+2H}|E^*(x)|^k\d x\cr&
\ll V^{-2k-2}L^{k+1}\int_{T-H}^{T+2H}|E^*(x)|^k\d x.
\cr}\leqno(3.9)
$$
This is somewhat sharper than the bound proved by the author in [6, Part II],
which contained $T^\e$ instead of a log-power, and the result was stated for the
``long'' interval $[T, 2T]$. The bound in (2.1) follows if the integral on the
left-hand side  is split into $O(\log T)$ subintegrals where $V \le |\zt|
\le 2V$. Denoting each such integral as $I_V$, we estimate it as
$$
I_V \ll \sum_{r=1}^{R_V}|\z(\hf+it_r)|^{2k+2} \ll R_VV^{2k+2}
\ll L^{k+1}\int_{T-H}^{T+2H}|E^*(x)|^k\d x,
$$
where the points $t_r$ are chosen in such a way that
$|t_r-t_s| \ge 1$ for $r\ne s$. Then (2.1) follows at once.

\medskip
To prove (2.2) we need ($C$ denotes generic positive constants)

\medskip
LEMMA 2. For $T^\e \le G= G(T)\le T, t \asymp T, L = \log T$ we have
$$
E^*(t) \le \frac{C}{G}\int_t^{t+G}\f_+(u) E^*(u)\d u + CGL,\leqno(3.10)
$$
and
$$
E^*(t) \ge \frac{C}{G}\int_{t-G}^{t}\f_-(u) E^*(u)\d u - CGL.\leqno(3.11)
$$
Here $\f_+$ is a non-negative, smooth function supported in $[t, t +G]$
such that $\f_+(u) = 1$ for $t+G/4 \le u \le t+3G/4$. Similarly, in (3.11) $\f_-$
is a non-negative, smooth function supported in $[t-G, t +G]$
such that $\f_-(u) = 1$ for $t-3G/4 \le u \le t-G/4$.
\medskip
The proof of these inequalities is similar, so it suffices only to prove
(3.10). From (1.2) we have, for $0 \le u\ll T$,
$$\eqalign{\cr
0 &\le \int_T^{T+u}|\zt|^2\d t = (T+u)\Bigl(\log\bigl({T+u\over2\pi}\bigr)
+ 2\gamma-1\Bigr) \cr&-  T\Bigl(\log\bigl({T\over2\pi}\bigr)
+ 2\gamma-1\Bigr) + E(T+u) - E(T).\cr}
$$
By the mean-value theorem this implies
$$
E(T) \le E(T+u) + O(u\log T),
$$
giving by integration and change of notation
$$
E(t) \le {C\over G}\int_t^{t+G} \f_+(u)E(u)\d u + CG\log T \quad(1 \ll G\ll T,\,C>0,\,t\asymp T).
\leqno(3.12)
$$
By using (3.2) again it is established that, for $T^\e \le G\le T, t \asymp T$,
$$
2\pi\D^*\Bigl({t\over2\pi}\Bigr)
= {C\over  G }\int_t^{t+G}\f_+(u)\D^*\Bigl({u\over2\pi}\Bigr)\d u+
O(G\log T). \leqno(3.13)
$$
Therefore by combining (3.12) and (3.13) one obtains (3.10), since
$$
E^*(t) = E(t) - 2\pi\D^*\Bigl({t\over2\pi}\Bigr).
$$

\medskip
In proving (2.2) we use (3.10) if $E^*(t) >0$, and (3.11) otherwise. Suppose
$E^*(t) >0$. Then by  integrations by parts we obtain from (3.10)
$$
\eqalign{
E^*(t) & \le \frac{C}{G}\int_t^{t+G}\f_+(u) E^*(u)\d u + CGL\cr&
= \frac{C}{G}\int_0^u E^*(v)\d v\cdot\f_+(u)\Bigg|_{u=t}^{t+G} + CGL
- \frac{C}{G}\int_t^{t+G}\f_+'(u)\int_0^u E^*(v)\d v\d u\cr&
= -\frac{C}{G}\int_t^{t+G}\left(\frac{3\pi}{4}u + R(u)\right)\f_+'(u)\d u+ CGL\cr&
= -\frac{C3\pi}{4H}\left(u\f_+(u)\Bigg|_{u=t}^{t+H}
- \int\limits_t^{t+G}\f_+(u)\d u\right) -\frac{C}{G}\int\limits_t^{t+G}R(u)\f_+'(u)\d u+ CGL
\cr&
= O(GL)  -\frac{C}{G}\int_t^{t+G}R(u)\f_+'(u)\d u.\cr}\leqno(3.14)
$$
Combining (3.14) with the corresponding lower bound and using the fact that
$$
\f_\pm'(u) \ll \frac{1}{G},
$$
 it follows that we have proved

\medskip
LEMMA 3. For $T^\e \le G = G(T)\le T,\, t\asymp T,$ we have
$$
|E^*(t)| \ll \frac{1}{G^2}\int_{t-G}^{t+G}|R(u)|\d u + CGL.\leqno(3.15)
$$

\medskip
If we suppose that $R(T) \ll_\e T^{\a+\e}$ then from (3.15), (3.5) of Lemma 1 and (3.1) we obtain
$$
\zt \ll_\e t^{\a/4+\e},  \quad E^*(T) \ll_\e T^{\a/2+\e},\leqno(3.16)
$$
so  that with the value $\a = 593/912 = 0.6502129...$ (see (1.11)) we have the bounds
$$
\eqalign{
\zt &\ll_\e |t|^{593/3648+\e},  \quad 593/3648 = 0,16255\ldots,\,\cr E^*(T) &\ll_\e T^{593/1824+\e},
\quad 593/1824 = 0,32510\ldots\,.\cr}\leqno(3.17)
$$
If  the conjectural $\a=\hf$ held ($\a < \hf$ is impossible by (1.14)),
then we would obtain from (3.15)
$$
\zt \ll_\e |t|^{1/8+\e},\qquad
E^*(T) \ll_\e T^{1/4+\e},
$$
which is out of reach by present day methods.
See (4.5) for the best known bound for $\zt$; the best known exponent for $E^*(T)$ is
$131/416 = 0.31490\ldots\,$. This exponent was proved for $E(T)$ by N. Watt [20], but
since the same exponent holds for $\D(x)$ and $\D^*(x)$, it holds for $E^*(T)$ as well.
Thus, although the bounds in (3.17) are non-trivial, they are not the best ones known
at present.
\medskip
We return now to our proof of (2.2).
Suppose now that $|E^*(t)| \ge V$ on a set of points $\{t_r\}_{r=1}^{\Cal R}$  lying in
$[T, T+H]$ and spaced at least $CG$ apart. We take $G = \delta V/L\;(< H)$ for sufficiently
small $\delta >0$. Then from (3.15) we have, for a representative set of the $t_r$'s such that
the intervals $(t_r-G,\,t_r+G)$ are disjoint,
$$
\eqalign{
{\Cal R}V^3L^{-2} &\ll \sum_{r=1}^{\Cal R}\int_{t_r-G}^{t_r+G}|R(u)|\d u\cr&
\ll \sum_{r=1}^{\Cal R}\left(\int_{t_r-G}^{t_r+G}|R(u)|^k\d u
\right)^{\frac{1}{k}}G^{1-\frac{1}{k}}\cr&
\ll \left(\sum_{r=1}^{\Cal R}\int_{t_r-G}^{t_r+G}|R(u)|^k\d u
\right)^{\frac{1}{k}}(RG)^{1-\frac{1}{k}},\cr}
$$
on applying H\"older's inequality for integrals. Since
the intervals $(t_r-G, t_r+G)$ are disjoint, and their union
is contained in $[T-H,\,T+2H]$, the preceding bound
gives us
$$
{\Cal R} \ll \int_{T-H}^{T+2H}|R(u)|^k\d u\cdot V^{-3k}L^{2k}G^{k-1},
$$
which simplifies to
$$
{\Cal R} \ll \int_{T-H}^{T+2H}|R(u)|^k\d u\cdot V^{-1-2k}L^{k+1}.\leqno(3.18)
$$
Splitting $\int\limits_T^{T+H}|E^*(t)|^{2k}\d t$ into $O(\log T)$ integrals $I_V$ where
$$
V \le |E^*(t)|\le 2V,
$$
we estimate each of these integrals by (3.18), keeping in mind that $V \le T^{1/3} \ll H$.
The bound in (2.2) follows at once.

\medskip
An obvious corollary of Theorem 1 is that
$$
 \int_{T}^{T+H}|\zt|^{4k+2}\d t \ll (\log T)^{3k+4} \int_{T-2H}^{T+4H}|R(t)|^k\d t
 \quad(T^{1/3} \ll H\ll T).\leqno(3.19)
$$
From (1.17) and (3.19) with $k=2$ we obtain
$$
\int_{T-H}^{T+H}|\zt|^{10}\d t \ll_\e T^\e(HT + T^{5/3})\quad(T^{1/3} \ll H\ll T).\leqno(3.20)
$$
It seems that this bound is new in the range when $H$ is close to $T^{1/3}$. It gives, by (3.1),
the classical bound $\zt \ll_\e |t|^{1/6+\e}$.

\bigskip
We shall now pass  to the proof of Theorem 2.
To obtain (2.4) we use (3.5) of Lemma 1 with $G \equiv H$.
This gives, for fixed $k\in\NN, T^\e \le H = H(T)\le T$,
$$
\eqalign{&
\int_T^{2T}\Bigl(
\int_{t-H}^{t+H}|\zu|^2\d u\Bigr)^k\d t \cr&
\ll H^{-k}\int_T^{2T}\left(\int_{-HL}^{HL}|E^*(t+x)|\E^{-(x/H)^2}\d x)\right)^k\d t +TH^kL^k.
\cr}\leqno(3.21)
$$
H\"older's inequality for integrals shows that the integral on the
right-hand side of (3.21) is
$$
\eqalign{&
\le \int_T^{2T}\int_{-HL}^{HL}|E^*(t+x)|^k\E^{-(x/H)^2}\d x\cdot\left(\int_{-HL}^{HL}
\E^{-(x/H)^2}\d x\right)^{k-1}\d t\cr&
\ll H^{k-1}\int_{-HL}^{HL}\E^{-(x/H)^2}\left(\int_{T-HL}^{2T+HL}|E^*(t+x)|^k\d t\right)\d x.\cr}
\leqno(3.22)
$$
From (3.21) and (3.22) we obtain (2.4) if we take into account (2.3). Note that the constant $A(k)$
in (2.3) actually must satisfy $A(k) \ge 1 + k/6$ for any $k \ge 1$, and not necessarily when
$k$ is an integer. If $k \ge 2$, then by H\"older's inequality for integrals
$$
\int_T^{2T}|E^*(t)|^2\d t \le \left(\int_T^{2T}|E^*(t)|^k\d t\right)^{2/k}T^{1-2/k},
$$
and the desired bound for $A(k)$ follows from the mean square formula (1.5). If $1\le k\le 2$
then it follows in a similar fashion from (1.5) and (1.7). We remark that if
$A(k) = 1 + k/6$ holds for some $k$, then (2.1) and (3.1) yield the bound
$$
\zt \;\ll_\e\; |t|^{\frac{k+6}{12(k+1)}+\e},
$$
and this improves the exponent $32/205 = 0,15609\ldots$ (see (4.5)) for $k\ge 5$, since for $k=5$ it gives
$11/72 = 0,152777\ldots$.

\bigskip
It remains to prove Theorem 3. We begin by noting that
 the author in [11] proved the following result, which improves on an
earlier result of M. Jutila [16]: If
$1 \ll U = U(T) \le \hf {\sqrt{T}}$, then we have ($c_3 = 8\pi^{-2}$)
$$\eqalign{
\int_T^{2T}\Bigl(\D(x+U)-\D(x)\Bigr)^2\d x & = TU\sum_{j=0}^3c_j\log^j
\Bigl({\sqrt{T}\over U}\Bigr) \cr&
+ O_\e(T^{1/2+\e}U^2) + O_\e(T^{1+\e}U^{1/2}),\cr}\leqno(3.23)
$$
a similar result being true
if $\D(x+U)-\D(x)$ is replaced by $E(x+U)-E(x)$, with different constants $c_j$ ($c_3>0$).
But the  integral in (2.7) can be reduced to the evaluation of the mean square
of $E(t+h)-E(t-h)$, since by (1.2) one has
$$
\int\limits_{t-H}^{t+H}|\zt|^2\d t = E(t+H)-E(t-H) + 2H\left(\log\bigl(\frac{t}{2\pi}\bigr) + 2\gamma\right)
+ O\left(\frac{H^2}{T}\right).\leqno(3.24)
$$
Therefore
$$
\int_T^{2T}\left(\int_{t-H}^{t+H}|\zu|^2\d u\right)^2\d t = I_1 + 2I_2 + I_3,
$$
say, where
$$
\eqalign{
I_1 &:= \int_T^{2T}\left(E(t+H)-E(t-H)\right)^2\d t =
\int_{T+H}^{2T+2H}\left(E(x+2H)-E(x)\right)^2\d x,\cr
I_2 &:= \int_T^{2T}2H\left(\log \frac{t}{2\pi} + 2\gamma+ O\left(\frac{H}{T}\right)\right)
(E(t+H)-E(t-H))\d t,\cr
I_3 &:= \int_T^{2T}4H^2\left(\log \frac{t}{2\pi} + 2\gamma+ O\left(\frac{H}{T}\right)\right)^2\d t.
\cr}\leqno(3.25)
$$
To evaluate $I_1$ we write
$$
I_1 = \int_{T+H}^{2T+2H} = \int_T^{2T} + \int_{2T}^{2T+2H} - \int_T^{T+H} = J_1 + J_2 - J_3,
$$
say. By trivial estimation, in view of $E(t) \ll t^{1/3}$ (see e.g., [4, Ch. 15]), it follows that
$$
J_2 - J_3 \;\ll\; HT^{2/3}.
$$
To evaluate $J_1$ we use the analogue of (3.23) (with $U=2H$) for $E(x+U)-E(x)$. This gives, with
suitable constants $d_j\;(d_3>0)$ and $1 \ll H \ll \sqrt{T}$,
$$
J_1 = TH\sum_{j=0}^3d_j\log^j\Bigl({\sqrt{T}\over 2H}\Bigr)
+ O_\e(T^{1/2+\e}H^2) + O_\e(T^{1+\e}H^{1/2}).
$$
One can evaluate $I_3$ in a straightforward way to obtain
$$
\eqalign{
I_3 &= 4H^2\int_T^{2T}\left(\log^2\bigl(\frac{t}{2\pi}\bigr) + 4\gamma^2 + 4\gamma
+ \log\bigl(\frac{t}{2\pi}\bigr) +O\left(\frac{H\log T}{T}\right)\right) \d t\cr&
= H^2T(4\log^2T + e_1\log T + e_0) + O(H^3\log T)\cr}
$$
with suitable constants $e_0$ and $e_1$.
\medskip
Finally to bound $I_2$ we invoke the result of J.L. Hafner and the author [2], namely
$$
E_1(T) := \int_2^T E(u)\d u = \pi T + O(G(T)), \quad G(T) = O(T^{3/4})\quad(T>2).\leqno(3.26)
$$
Actually in [2] an explicit expression is given for $G(T)$ (from which one can deduce that
$G(T) = \Omega_\pm(T^{3/4})$). Thus from (3.25), (3.26) we obtain, on integrating by parts,
$$
\eqalign{
I_2&= 2H\Biggl\{\Bigl(E_1(t+H)-E_1(t-H)\Bigr)\left(\log \frac{t}{2\pi}
+ 2\gamma\right)\Biggr\}\Biggl|_{t=T}^{2T}
\cr&
- 2H\int_T^{2T}(E_1(t+H)-E_1(t-H))\frac{\d t}{t} + O(H^2T^{1/3})\cr&
= O(H^2\log T) + O(HT^{3/4}\log T) +  O(H^2T^{1/3}) = O(HT^{3/4}\log T)\cr}
$$
in view of the range for $H$, namely $T^\e \le H = H(T) \le T^{1/2-\e}$.

\medskip
Combining the expressions for $I_1, I_2$ and $I_3$ we obtain (2.7), which in the range
$T^\e \le H = H(T) \le T^{1/2-\e}$ provides an asymptotic formula for the integral in question.
Note that in this range $HT^{3/4}L\ll T^{1+\e}H^{1/2}$, so only the error
terms in (2.7) remain. For $T^{1/2-\e} \le H \le T$ the upper bound in (2.6) follows easily from
(2.4) and $A(2) \le 4/3$, see (4.1).

\medskip
It remains yet to prove (2.5). Note that, by (3.24), the integral in question
is easily seen to be equal to
$$
2H\int_T^{2T}\Bigl(\log \frac{t}{2\pi}+ O\Bigl(\frac{H}{T}\Bigr)\Bigr)\d t +
\int_T^{2T}\Bigl(E(t+H)-E(t-H)\Bigr)\d t.
\leqno(3.27)
$$
But by using (3.26) again it is seen that (3.27) reduces to
$$
2H\left(T\log \Bigl(\frac{2T}{\E\pi} \Bigr)\right) + O(H^2) + O(T^{3/4}).
$$
Hence, for $T^\e \le H =H(T)\le T$,
$$
\int_T^{2T}\Bigl(\int_{t-H}^{t+H}|\zu|^2\d u\Bigr)\d t =
2H\left(T\log \Bigl(\frac{4T}{\E} \Bigr)\right) + O(H^2) + O(T^{3/4}),
$$
as asserted by (2.5).

\bigskip
\head
3. Some corollaries and remarks
\endhead

\bigskip
If $A(k)$ is defined by (2.3), then from (1.7)--(1.9)  we have
$$
A(2) \le \frac{4}{3},\quad A(3) \le \frac{3}{2},\quad A(4) \le \frac{7}{4},\quad A(5) \le 2.
\leqno(4.1)
$$
We also have $A(1) \le 7/6$ by $A(2) \le 4/3$ and the Cauchy-Schwarz inequality.
Then (with $H=T$) (2.1) of Theorem 1 yields
$$
\eqalign{&
\int_0^T|\zt|^8\d t \ll_\e T^{3/2+\e},\cr& \int_0^T|\zt|^{10}\d t \ll_\e T^{7/4+\e},\cr&
\int_0^T|\zt|^{12}\d t \ll_\e T^{2+\e},\cr}\leqno(4.2)
$$
with $k=3,4,5$, respectively.
The bounds in (4.2) (up to $T^\e$, which can be replaced by a log-factor) are the
sharpest known bounds for the moments in question (see e.g., [4, Chapter 8]).

On the other hand, by using (1.4), we also have from (2.1)
$$
\int_T^{T+H}|\zt|^6\d t \ll_\e HT^{1/3}\log^{12}T + T^{1+\e}\quad(T^{1/3} \le H \le T).\leqno(4.3)
$$
Although this is not trivial, it can be improved if one uses the bound of H. Iwaniec [13]
$$
\int_T^{T+H}|\zt|^4\d t \ll_\e T^\e(H + TH^{-1/2})\quad(T^\e \le H \le T). \leqno(4.4)
$$
The bound in (4.4) was obtained by sophisticated methods from the spectral theory
of the non-Euclidean Laplacian, and if coupled with the best known bound of
M.N. Huxley [3] for $|\zt|$, namely
$$
\zt \ll_\e |t|^{32/205 + \e}, \quad 32/205 = 0,15609\ldots\,, \leqno(4.5)
$$
one gets an improvement of (4.3). Note that the famous, yet unsettled Lindel\"of conjecture
states that, instead of (4.5), one has $\zt \ll_\e |t|^\e$.

\medskip
If we combine (1.16) and (2.2) (with $k=2$), it follows that
$$
\int_T^{T+H}|E^*(t)|^4\d t \ll_\e HT\log^7T + T^{5/3+\e}\quad(T^{1/3+\e}\le H \le T).\leqno(4.6)
$$
The bound in (4.6) does not follow from (1.9), as it is better for $T^{2/3}\le H\le T^{3/4}$.

\medskip
As a corollary to Theorem 2, we obtain with (4.1)
$$
\eqalign{&
\int_T^{2T}\left(\int_{t-H}^{t+H}|\zt|^2\d u\right)^3\d t \ll_\e T^{3/2+\e} + TH^3L^3,\cr&
\int_T^{2T}\left(\int_{t-H}^{t+H}|\zt|^2\d u\right)^4\d t \ll_\e T^{7/4+\e}+ TH^4L^4,\cr&
\int_T^{2T}\left(\int_{t-H}^{t+H}|\zt|^2\d u\right)^5\d t \ll_\e T^{2+\e} + TH^5L^{5}.
\cr}
\leqno(4.7)
$$
All the bounds in (4.7) are valid for $T^\e \le H\le T$, but as we have (see e.g., K. Ramachandra [18])
$$
\int_{t-H}^{t+H}|\zt|^{2k}\d t \;\gg_k\; H(\log H)^{k^2}\qquad(\log\log T\ll H \le T,\,k\in\NN),
$$
we have the expected upper bounds $T(HL)^{m}\;(m = 3,4,5)$ for the integrals in (4.7). Indeed, we
obtain from (4.7)
$$
\eqalign{&
\int_T^{2T}\left(\int_{t-H}^{t+H}|\zu|^2\d u\right)^3\d t \;\ll\; TH^3L^3\quad(H \ge T^{1/6+\e}),\cr&
\int_T^{2T}\left(\int_{t-H}^{t+H}|\zu|^2\d u\right)^4\d t \;\ll\;  TH^4L^4\quad(H \ge T^{3/16+\e}),\cr&
\int_T^{2T}\left(\int_{t-H}^{t+H}|\zu|^2\d u\right)^5\d t \;\ll\; TH^5L^{5}\quad(H \ge T^{1/5+\e}).
\cr}
\leqno(4.8)
$$
The bounds in (4.8) seem to be the best unconditional bounds yet.

\medskip
Note that for the analogous, but less difficult, problem of moments of
$$
J_k(t,G) := {1\over\sqrt{\pi}G}
\int_{-\infty}^\infty |\z(\hf + it + iu)|^{2k}{\roman e}^{-(u/G)^2}\d u
\qquad(t \asymp T, T^\e \le G \ll T),$$
where $k$ is a natural number, we refer the reader to the author's work [7]. Not only do we have
$$
\int_{T-G}^{T+G}|\zt|^{2k}\d t= \int_{-G}^G|\z(\hf + iT + iu)|^{2k}\d u
\le \sqrt{\pi}{\roman e}G\,J_k(T,G),
$$
but the presence of the smooth Gaussian exponential factor in $J_k(T,G)$ facilitates
the ensuing estimations. We have (this is [7, Theorem 1])
$$
\int_T^{2T}J_1^m(t,G)\d t \ll_\e T^{1+\e}\leqno(4.9)
$$
for $T^\e \le G \le T$ if $m = 1,2$; for $T^{1/7+\e} \le G \le T$ if
$m=3$, and for $T^{1/5+\e} \le G \le T$ if $m=4$; and these bounds were sharpened in [10]
to $T^{7/36} \le G \le T$ when $m =4$, $T^{1/5} \le G \le T$ when $m =5$ and
$T^{2/9} \le G \le T$ when $m=6$. The bounds in (4.9) can be
compared to those in (4.8).

\bigskip
We remark that in [8] the author proved that
$$
\int_T^{2T}\left(\int_{t-G}^{t+G}|\zu|^4\d u\right)^2\d t \ll_\e G^2T^{1+\e}\leqno(4.10)
$$
for $T^{1/2} \le G = G(T) \ll T$. In fact, (4.10) is connected with the following, more general
result (Theorem 1 of [8]):
Let $T < t_1 < t _2 < \ldots < t_R <
2T$, $t_{r+1} - t_r \ge G$ for $r = 1, \ldots\,, R-1$. If, for
fixed $m,k\in \NN$, we have
$$
\int_T^{2T}\left({1\over G}\int_{t-G}^{t+G}|\zeta(\hf+iu)|^{2k}\d
u\right)^m\d t \;\ll_\e\; T^{1+\e} \leqno(4.11)
$$
for $T^{\a_{k,m}} \le G = G(T) \ll T$ and $0\le \a_{k,m}\le 1$,
then
$$
\sum_{r=1}^R\int_{t_r-G}^{t_r+G}|\zt|^{2k}\d t \ll_\e
(RG)^{m-1\over m} T^{{1\over m}+\e}.
$$
In this notation, (4.10) is implied by $\a_{2,2} = \hf$. In fact, if (4.11) holds, then we have
$$
\int_0^T|\zt|^{2km}\d t \;\ll_\e\;
T^{1+(m-1)\a_{k,m}+\e}.
$$
Non-trivial bounds of the type (4.11) (with $0\le \a_{k,m}\le 1$) are hard to obtain when $m>2$ or $k>2$.

\vfill
\eject
\Refs

\item{[1]} F.V. Atkinson, {\it The mean value of the Riemann zeta-function},
Acta Math. {\bf81}(1949), 353-376.

\item{[2]} J.L. Hafner and A. Ivi\'c, {\it On the mean square of the Riemann
zeta-function on the critical line}, J. Number Theory
   {\bf 32}(1989), 151-191.
\item{[3]} M.N. Huxley,
{\it Exponential sums and the Riemann zeta function V},
Proc. London Math. Soc. (3) {\bf 90}(2005), 1-41.

\item{[4]} A. Ivi\'c, {\it The Riemann zeta-function}, John Wiley \&
Sons, New York, 1985 (2nd ed. Dover, Mineola, New York, 2003).

\item {[5]} A. Ivi\'c,  {\it Mean values of the Riemann zeta-function},
LN's {\bf 82},  Tata Inst. of Fundamental Research,
Bombay,  1991 (distr. by Springer Verlag, Berlin etc.).

\item{[6]} A. Ivi\'c, {\it On the Riemann zeta-function and the divisor problem},
Central European J. Math. {\bf(2)(4)}\ (2004), 1-15,  II, ibid.
{\bf(3)(2)}(2005), 203-214,  III, Annales Univ.
Sci. Budapest, Sect. Comp. {\bf29}(2008), 3-23,
and IV, Uniform Distribution Theory {\bf1}(2006), 125-135.

\item{[7]} A. Ivi\'c, {\it On moments of $|\zeta({1\over2}+it)|$ in short intervals},
Ramanujan Math. Soc. LNS{\bf2},
The Riemann zeta function and related themes: Papers in honour of Professor
K. Ramachandra (Proc. Conference held at Bangalore 13-15 Dec. 2003, eds.
R. Balasubramanian and K. Srinivas), 2006, 81-97.

\item{[8]} A. Ivi\'c, {\it On sums of integrals of powers of the zeta-function
in short intervals},
Multiple Dirichlet Series, Automorphic Forms, and Analytic Number Theory
(eds. S. Friedberg et al.), Proc. Symposia Pure Math.
Vol. {\bf75}, AMS, Providence, Rhode Island, 2006, pp. 231-242.

\item{[9]} A. Ivi\'c, {\it On the mean square of the zeta-function and
the divisor problem}, Annales  Acad. Scien. Fennicae Mathematica {\bf23}(2007), 1-9.

\item{[10]} A. Ivi\'c, {\it Some remarks on the moments of $|\zeta({1\over2}+it)|$
in short intervals}, Acta Math. Hung. {\bf119}(2008), 15-24.

\item{[11]} A. Ivi\'c, {\it On the divisor function and the Riemann zeta-function
in short intervals},
The Ramanujan Journal {\bf19}(2)(2009), 207-224.

\item{[12]} A. Ivi\'c, {\it On some mean square  estimates for the
zeta function in short intervals},
 to appear in Annales Univ.
Sci. Budapest, Sect. Comp., preprint in arXiv:1212.0660.

\item{[13]} H. Iwaniec, {\it Fourier coefficients of Cusp Forms and
the Riemann Zeta-Function}, Expos\'e No. {\bf18},
S\'eminaire de Th\'eorie des Nombres, Universit\'e Bordeaux, 1979/80.

\item{[14]} M. Jutila,{\it Riemann's zeta-function and the divisor problem},
Arkiv Mat. {\bf21}(1983), 75-96 and II, ibid. {\bf31}(1993), 61-70.

\item{[15]} M. Jutila, {\it On a formula of Atkinson},
Topics in classical number theory, Colloq. Budapest 1981,
Vol. I, Colloq. Math. Soc. J\'anos Bolyai {\bf34}(1984), 807-823.

\item{[16]} M. Jutila, {\it On the divisor problem for short intervals},
Ann. Univer. Turkuensis Ser. {\bf A}I {\bf186}\break(1984), 23-30.

\item{[17]} T. Meurman, {\it A generalization of Atkinson's formula to
$L$-functions}, Acta Arith. {\bf47}(1986), 351-370.

\item{[18]} K.  Ramachandra, {\it On the mean-value and omega-theorems
for the Riemann zeta-function}, LN's {\bf85}, Tata Inst. of Fundamental Research
(distr. by Springer Verlag, Berlin etc.), Bombay, 1995.

\item{[19]} P. Shiu, {\it A Brun-Titchmarsh theorem for multiplicative functions},
J. Reine Angew. Math. {\bf313}(1980), 161-170.

\item{[20]} N. Watt, {\it A note on the mean square of $|\zeta (\frac12 + it)|$},
J. Lond. Math. Soc., II. Ser. {\bf82}, No. 2, (2010), 279-294.
\vfill

\endRefs

\enddocument

\bye